\documentclass[12pt]{article}

\usepackage{amstext    }
\usepackage{amsthm    }
\usepackage{a4}
\usepackage[mathscr]{eucal}
\usepackage{mathrsfs}

\usepackage{amsmath}
\usepackage{amssymb}
\usepackage{amscd}

\newtheorem{theorem}{Theorem}[section]
\newtheorem{definition}[theorem]{Definition}
\newtheorem{proposition}[theorem]{Proposition}
\newtheorem{corollary}[theorem]{Corollary}
\newtheorem{lemma}[theorem]{Lemma}

\newtheorem{example}[theorem]{Example}

\newcommand{\cali}[1]{\mathscr{#1}}

\newcommand{\Ec}{\cali{E}}

\newcommand{\C}{\mathbb{C}}

\renewcommand{\P}{\mathbb{P}}

\title{Analytic multiplicative cocycles over holomorphic
  dynamical systems}
\author{Tien-Cuong Dinh}

\begin{document}
\maketitle

\begin{abstract}
We prove some properties of analytic
multiplicative and sub-multiplicative 
cocycles. The results allow to
construct natural invariant analytic sets associated to complex
dynamical systems.    
\end{abstract}

\noindent {\bf AMS classification:} 37F, 32H

\noindent {\bf Key-words:} multiplicative cocycle, sub-multiplicative
cocycle, exceptional set


\section{Introduction} \label{section_introduction}

In this paper we study some sub-multiplicative cocycles over complex dynamical systems.
The notion of sub-multiplicative cocycle is classical. It lies at the root of the
definition of Lyapunov exponents. In \cite{Favre} cocycles which are
upper semi-continuous with respect to the analytic Zariski topology, were introduced and
studied. Applications to holomorphic dynamics were subsequently given
in  \cite{Favre2, FavreJonsson}.
Our purpose here is to introduce a new approach closer in spirit to \cite{DinhSibony1,
  DinhSibony3}. The properties
we obtain strengthen former results of Favre. In Section \ref{section_dynamics} below, we give an application
to equidistribution problems of our main theorem, see also \cite{DinhSibony1,
  DinhSibony3}.

Let $X$ be an irreducible compact complex space of dimension
$k$, not necessarily
smooth. Let $f:X\rightarrow X$ be an open holomorphic map. 
Consider the dynamical system associated to $f$ and define 
$f^n:=f\circ\cdots\circ f$, $n$ times, the iterate of order $n$ of $f$.

\begin{definition} \label{def_cocycle}
\rm
A sequence of functions $\kappa_n:X\rightarrow [1,+\infty[$,
$n\geq 0$, is an {\it analytic sub-multiplicative} (resp. {\it
multiplicative}) {\it cocycle} if for any $n,m\geq 0$ and
any $x\in X$, we have
\begin{enumerate}
\item $\kappa_n$ is upper semi-continuous for the Zariski
topology on $X$ and
 $\min_X \kappa_n=1$;
\item  $\kappa_{n+m}(x)\leq \kappa_n(x)\kappa_m(f^n(x))$
(resp. $\kappa_{n+m}(x)= \kappa_n(x)\kappa_m(f^n(x))$).
\end{enumerate}
\end{definition}

The first property implies that $\kappa_n$ is bounded
from above and equal to 1 out of a finite or countable union of
proper analytic subsets of $X$. Moreover, for any $\delta>1$, 
$\{\kappa_n\geq \delta\}$ is a proper
analytic subset of $X$ and there is $\delta'<\delta$ such that
$\{\kappa_n<\delta\}$ is contained in 
$\{\kappa_n<\delta'\}$. The second property defines classical
multiplicative and sub-multiplicative cocycles, see \cite{Walters}. It 
implies that if $M$ is a constant such that $\kappa_1\leq
M$, then $\kappa_n\leq M^n$ for $n\geq 1$. 

The new tool we introduce is to define for any analytic
sub-multiplicative  cocycle $(\kappa_n)_{n\geq 0}$, its extension to the negative integers:
$$\kappa_{-n}(x):=\max_{y\in f^{-n}(x)}\kappa_n(y).$$
Since $\{\kappa_{-n}\geq \delta\}=f^n\{\kappa_n\geq\delta\}$, the function
$\kappa_{-n}$ is upper semi-continuous in the
Zariski sense. Here is our main result.

\begin{theorem} \label{th_cocycle}
The sequence of functions $[\kappa_{-n}]^{1/n}$ converges
pointwise to a function $\kappa_-$. Moreover, for every $\delta>1$,
 the level set $\{\kappa_-\geq\delta\}$ is a proper
analytic subset of $X$ which is invariant under $f$ and is contained in the orbit of
$\{\kappa_n\geq \delta^n\}$ for every $n\geq 0$. In particular, $\kappa_-$ is upper
semi-continuous in the Zariski sense.
\end{theorem}

Using the function $\kappa_-$, we obtain a more uniform version of
\cite[Th. 2.3.5]{Favre}
for sub-multiplicative cocycles, which is proved to be useful in
applications, see Section
\ref{section_dynamics} below.

\begin{corollary} \label{cor_cocycle}
If $\limsup [\kappa_n(x)]^{1/n}\geq\delta$, then $x$ is sent by some
  $f^n$ to a point in the invariant
analytic set $\{\kappa_-\geq \delta\}$. If $(\kappa_n)$ is
multiplicative, then $[\kappa_n]^{1/n}$ converge pointwise to a
function $\kappa_+$ such that $\kappa_+\circ f=\kappa_+$.
\end{corollary}

The new and useful property here is that the analytic set
$\{\kappa_-\geq \delta\}$ does not depend on $x$.
Note that
one can extend the study to the case
where $f:X\rightarrow X$ is a
continuous open map on an irreducible Notherian topological space of
finite dimension, see \cite{Favre,Hartshorne}. In particular, the above results
hold for open real analytic maps and for open algebraic maps on
quasi-projective varieties. We will discuss the case of dominant 
meromorphic maps on K\"ahler compact
manifolds in Section \ref{section_dynamics}.

Note also that the condition $\min\kappa_n=1$ in Definition
\ref{def_cocycle} can be replaced by $\kappa_n\geq c^n$ for some
constant $c>0$. Indeed, let $Y$ be a compact analytic set and $b$ a
point in $Y$. Consider the
map $f':X\times Y\rightarrow X\times Y$ defined by $f'(x,y):=(f(x),y)$ and
the cocycle $\kappa'_n(x,y):=1$ if $y\not=b$
and $\kappa'_n(x,b):=2^nc^{-n}\kappa_n(x)$. Theorem \ref{th_cocycle} and
Corollary \ref{cor_cocycle} can be
applied to $f'$ and $\kappa'_n$. We deduce analogous results for 
$\delta>\delta^*:=\min\kappa_-$ (in Theorem \ref{th_cocycle} and
Corollary \ref{cor_cocycle} we have $\delta^*=1$).

\section{Asymptotic behavior of cocycles}

In this section, we prove the results stated
above. Observe that 
$(\kappa_{nl})_{n\geq 0}$ is an analytic (sub-)multiplicative cocycle for
$f^l$ with $l\geq 1$. We have the
following elementary lemmas. 

\begin{lemma} \label{lemma_cocycle_bis}
Let $l\geq 1$ be an integer. Then 
$$\limsup_{n\rightarrow\infty}[\kappa_{nl}]^{1/(nl)}=\limsup_{n\rightarrow\infty}[\kappa_n]^{1/n} \quad\mbox{and}\quad
\liminf_{n\rightarrow\infty}[\kappa_{nl}]^{1/(nl)}=\liminf_{n\rightarrow\infty}[\kappa_n]^{1/n}.$$
\end{lemma}
\proof 
If $M$ is a constant such that $\kappa_s\leq M$ for any
$0\leq s\leq l$, then
\begin{equation*}
\kappa_{m+s}(x)\leq \kappa_m(x)\kappa_s(f^m(x)) \leq M
\kappa_m(x). 
\end{equation*}
We deduce from these inequalities that
\begin{equation*}
 M^{-1}\kappa_{(n+1)l}\leq \kappa_{nl+s}\leq M \kappa_{nl}.
\end{equation*}
Hence
$$\limsup_{n\rightarrow\infty}[\kappa_{nl+s}]^{1/(nl+s)}=\limsup_{n\rightarrow\infty}[\kappa_{nl}]^{1/(nl)}$$
for every $0\leq s\leq l-1$. The first identity in the lemma follows. The last one 
is obtained in the same way.
\endproof

\begin{lemma}\label{lemma_cocycle_3}
Let $l\geq 1$ be an integer. Then 
$$\limsup_{n\rightarrow\infty}[\kappa_{-nl}]^{1/(nl)}=\limsup_{n\rightarrow\infty}[\kappa_{-n}]^{1/n} \quad\mbox{and}\quad
\liminf_{n\rightarrow\infty}[\kappa_{-nl}]^{1/(nl)}=\liminf_{n\rightarrow\infty}[\kappa_{-n}]^{1/n}.$$
\end{lemma}
\proof
As in Lemma \ref{lemma_cocycle_bis}, we get 
\begin{equation}
\kappa_{m+s}(x)\leq \kappa_s(x)\kappa_m(f^s(x))\leq  M
\kappa_m(f^s(x)). \label{eq_cocycle}
\end{equation}
Therefore, if $a\in X$ and $b\in f^{-m-s}(a)$ are such that
$\kappa_{m+s}(b)=\kappa_{-m-s}(a)$, then 
$$\kappa_{-m-s}(a)=\kappa_{m+s}(b)\leq M\kappa_m(f^s(b))\leq
M\max_{x\in f^{-m}(a)} \kappa_m(x)=M\kappa_{-m}(a).$$
We deduce that 
\begin{equation*}
 M^{-1}\kappa_{-(n+1)l}\leq \kappa_{-nl-s}\leq
M \kappa_{-nl}. 
\end{equation*}
The lemma follows.
\endproof

For any irreducible analytic subset $Y$ of $X$, define
$$\kappa_n(Y):=\min_{x\in Y} \kappa_n(x) \quad \mbox{and} \quad \kappa_{-n}(Y):=\min_{x\in Y}
\kappa_{-n}(x).$$ 
Since $\kappa_n$ and $\kappa_{-n}$ are upper
semi-continuous for the Zariski topology, the previous minimums always
exist. Moreover, we have $\kappa_n(Y)=\kappa_n(x)$ and
$\kappa_{-n}(Y)=\kappa_{-n}(x)$ for $x\in Y$ outside a finite or
countable union of proper analytic subsets of $Y$.

\begin{lemma} \label{lemma_cocycle}
Let $Y$ be a periodic irreducible analytic subset of $X$. Then
$[\kappa_n(Y)]^{1/n}$ converge to a constant $\kappa_+(Y)$.
Moreover,  for every $m\geq
0$, we have $\kappa_+(f^m(Y))=\kappa_+(Y)$ and $\liminf[\kappa_{-n}(Y)]^{1/n}\geq\kappa_+(Y)$. If
$\kappa_+(Y)\geq\delta$ then  for every $m\geq 0$,  $Y$ is contained in the orbit of
$\{\kappa_m\geq\delta^m\}$.
\end{lemma}
\proof Assume that $f^l(Y)=Y$ with $l\geq 1$. Using properties of
cocycles for a generic point in $Y$,
we have
$$\kappa_{(m+n)l}(Y)\leq \kappa_{ml}(Y)
\kappa_{nl}(f^{ml}(Y))=\kappa_{ml}(Y)\kappa_{nl}(Y).$$ 
Hence, the
sequence of $[\kappa_{nl}(Y)]^{1/(nl)}$ is decreasing and
converges to some constant $\kappa_+(Y)$ when $n$ tends to infinity. The first assertion
follows from Lemma \ref{lemma_cocycle_bis} applied to a generic
point in $Y$. 
By definition of $\kappa_{-n}$,  since $f^l(Y)=Y$, we have $\kappa_{-nl}(Y)\geq
\kappa_{nl}(Y)$.
This and Lemma \ref{lemma_cocycle_3} imply the assertion on $\kappa_{-n}$.

For the second assertion, applying (\ref{eq_cocycle}) for a
generic point in $Y$ or $f(Y)$ yields
$$M^{-1}\kappa_{n+1}(Y)\leq \kappa_n(f(Y))\leq M\kappa_{n-l+1}(f^l(Y))=
M\kappa_{n-l+1}(Y).$$ This and the first assertion imply that
$[\kappa_n(f(Y))]^{1/n}$ converge to $\kappa_+(Y)$. It follows
that $\kappa_+(f(Y))=\kappa_+(Y)$ and then
$\kappa_+(f^m(Y))=\kappa_+(Y)$ for every $m\geq 0$.

Assume that $\kappa_+(Y)\geq \delta$ and that $Y$ is not contained in 
the orbit of $\{\kappa_m\geq\delta^m\}$ for a fixed $m\geq 0$. Since $Y$
is periodic,
$f^i(Y)$ is not contained in $\{\kappa_m\geq\delta^m\}$ for every
$i\geq 0$. We then deduce that $\kappa_m(f^i(Y))\leq \delta_0^m$
for some constant $\delta_0<\delta$. This implies
$$\kappa_{nm}(Y)\leq \prod_{i=0}^{n-1}\kappa_m(f^{im}(Y))
\leq\delta_0^{nm}$$ that contradicts the inequality
$\kappa_+(Y)\geq\delta$.
\endproof

We have the following proposition.

\begin{proposition} \label{prop_cocycle}
Let $q$ be an integer with $0\leq q\leq k-1$. Let $\Omega$ be a dense
Zariski open subset of $X$ such that $f^{-1}(\Omega)\subset\Omega$.
Assume 
that $\{\kappa_{-n_0}\geq
\delta^{n_0}\}\cap\Omega$ is of dimension at most equal to $q$ for
some $n_0\geq 1$ and some $\delta>1$. Then there is an integer
$n_1\geq 1$ and 
an invariant proper 
analytic subset $\Sigma$ of $X$, possibly reducible, which is empty or of pure dimension $q$, 
such that
\begin{enumerate}
\item  
$\Sigma$ is contained in the orbit of $\{\kappa_n\geq
\delta^n\}$ for every $n\geq 0$;
\item 
$\liminf[\kappa_n(x)]^{1/n}\geq \delta$ and
$\liminf[\kappa_{-n}(x)]^{1/n}\geq \delta$ for $x\in\Sigma$;
\item $\{\kappa_{-n_1}\geq
\delta^{n_1}\}\cap\Omega\setminus\Sigma$ is of dimension $\leq q-1$ if
$q\geq 1$ and 
is empty if $q=0$.
\end{enumerate}
\end{proposition}
\proof
Let $\widetilde \Sigma$ denote the union of the irreducible components of dimension
$q$ in the closure of $\{\kappa_{-n_0}\geq \delta^{n_0}\}\cap\Omega$.
Recall that $\kappa_{-n_0}$ is upper semi-continuous and
that $\{\kappa_{-n_0}\geq \delta^{n_0}\}\cap\Omega$ is of dimension
less than or equal to $q$. Hence, there is a real number $1\leq\delta_0<\delta$ satisfying
$\kappa_{n_0}(V)\leq \delta_0^{n_0}$ for any irreductible
analytic subset $V$ of dimension $\geq q$ of $X$ such that $f^{n_0}(V)$ intersects
$\Omega\setminus\widetilde\Sigma$. We will cover $\widetilde \Sigma$
by three analytic sets $\Sigma$, $\Sigma'$ and $\Sigma''$.

Define $\Sigma$ as the union of the orbits of periodic irreducible
components $Y$ of $\widetilde \Sigma$ such that $\kappa_+(Y)\geq
\delta$. So, $\Sigma$ is invariant, i.e. $f(\Sigma)=\Sigma$. 
Note that $\Sigma$ is not necessarily contained in
$\widetilde\Sigma$. 
By Lemma \ref{lemma_cocycle}, we also have
$\kappa_+(f^i(Y))\geq\delta$ for every $i\geq 0$ and
$\Sigma$ satisfies the first two properties in the proposition. Let $\Sigma'$ be the union of
the orbits of periodic irreducible components $Y$ of
$\widetilde\Sigma$ with $\kappa_+(Y)<\delta$. Then, $\Sigma'$ is invariant
and by Lemma \ref{lemma_cocycle}, $\kappa_+(Z)<\delta$ for any
irreducible component $Z$ of $\Sigma'$. Therefore, we can find $m\geq 0$ and
$1\leq\delta_1<\delta$ such that $\kappa_{mn_0}(Z)\leq
\delta_1^{mn_0}$ for any irreducible component $Z$ of
$\Sigma'$. Finally denote by $\Sigma''$ the
union of the non-periodic irreducible components of
$\widetilde\Sigma$ and $l$ the number of  irreducible components of $\Sigma''$.

Choose an integer $N$ large enough so that 
$M^{l+m}\max(\delta_0,\delta_1)^{Nmn_0}<\delta^{Nmn_0}$ where $M$ denotes the maximal
value of $\kappa_{n_0}$ on $X$. It is enough to check the last property in the
proposition for $n_1:=Nmn_0$. For this purpose, we only have to show
that $\kappa_{-n_1}(V_0)<\delta^{n_1}$ for any irreducible analytic
set $V_0$ of dimension $q'\geq q$ which intersects $\Omega\setminus
\Sigma$.
 Consider such an analytic set $V_0$ and an arbitrary  sequence $V_{-Nm}$, $V_{-Nm+1}$, $\ldots$, $V_0$ of
irreducible analytic subsets of dimension $q'$ of $X$ such that
$V_{-i+1}=f^{n_0}(V_{-i})$.
Since $f^{-1}(\Omega\setminus\Sigma)\subset
\Omega\setminus\Sigma$ and since $f$ is surjective, all these analytic sets intersect $\Omega\setminus\Sigma$.
It is also clear that the considered sequence  contains 
at most $l$ elements which are  irreducible components of
$\Sigma''$. We show that $\kappa_{Nmn_0}(V_{-Nm})<\delta^{Nmn_0}$
which implies that $\kappa_{-n_1}(V_0)<\delta^{n_1}$ and completes the
proof. We distinguish three cases.

If $V_0$ is not an  irreducible  component of $\Sigma'$ (in particular
when $q'>q$), the sequence of $V_{-n}$
does not contain any  irreducible component of $\Sigma'$. It contains at most
$l$ elements which are  irreducible components of $\Sigma''$. 
The other components intersect $\Omega\setminus\widetilde\Sigma$. Therefore,
we have
$$\kappa_{Nmn_0}(V_{-Nm})\leq\prod_{n=0}^{Nm-1}
\kappa_{n_0}(V_{-Nm+n})\leq M^l\delta_0^{Nmn_0}<\delta^{Nmn_0}.$$
From now on, we only have to consider the case where $q'=q$.

If $V_{-Nm}$ is an  irreducible component of $\Sigma'$ then $V_{-i}$
is an  irreducible component
of $\Sigma'$ for $0\leq i\leq Nm$ since $\Sigma'$ is invariant. It
follows from properties of cocycles that
$$\kappa_{Nmn_0}(V_{-Nm}) \leq \prod_{i=0}^{N-1}
\kappa_{mn_0}(V_{-(N-i)m})\leq \delta_1^{Nmn_0}<\delta^{Nmn_0}.$$

Otherwise, let $0\leq s\leq N-1$ be the largest integer such that
$V_{-im}$ is an  irreducible  component of $\Sigma'$ for $i\leq s$. Then $V_{-n}$ is
not an  irreducible  component of $\Sigma'$ for $n\geq (s+1)m$. Therefore, since the
considered sequence contains at most $l$ elements which are
 irreducible components of $\Sigma''$, we have
\begin{eqnarray*}
\kappa_{Nmn_0}(V_{-Nm}) & \leq &  \prod_{n=0}^{Nm-sm-1}
\kappa_{n_0}(V_{-Nm+n})\prod_{i=0}^{s-1}\kappa_{mn_0}(V_{-{(s-i)m}})\\
& \leq & \prod_{n=0}^{Nm-(s+1)m-1}
\kappa_{n_0}(V_{-Nm+n})\prod_{Nm-(s+1)m}^{Nm-sm-1}
\kappa_{n_0}(V_{-Nm+n})\prod_{i=0}^{s-1}\kappa_{mn_0}(V_{-{(s-i)m}})\\
& \leq &  M^l\delta_0^{[Nm-(s+1)m]n_0}M^m\delta_1^{smn_0}<\delta^{Nmn_0}.
\end{eqnarray*}
This completes the proof.
\endproof

\noindent 
{\bf Proof of Theorem \ref{th_cocycle}.} Fix a real
number $\delta>1$. We apply Proposition \ref{prop_cocycle} to
$n_0:=1$, $q:=k-1$ and $\Omega=\Omega_0:=X$. Since $\min\kappa_1=1$, we can find an
invariant hypersurface $\Sigma_1$ and an integer $n_1$ such that
$\{\kappa_{-n_1}\geq\delta^{n_1}\}\setminus\Sigma_1$ is of dimension
at most equal to $k-2$. Moreover, $\liminf [\kappa_n]^{1/n}\geq \delta$ and
$\liminf [\kappa_{-n}]^{1/n}\geq\delta$ on $\Sigma_1$. The Zariski
open set $\Omega_1:=X\setminus \Sigma_1$ satisfies
$f^{-1}(\Omega_1)\subset \Omega_1$. Hence, we can apply
Proposition \ref{prop_cocycle} inductively on $q$. We find
invariant analytic subsets $\Sigma_{k-q}$ of pure dimension $q$,
$0\leq q\leq k-1$, such that $\liminf [\kappa_n]^{1/n}\geq \delta$
and $\liminf [\kappa_{-n}]^{1/n}\geq\delta$ on $\Sigma_q$. Moreover,
if $\Omega_{k-q}:=X\setminus \Sigma_1\cup\ldots\cup \Sigma_{k-q}$, then
$\{\kappa_{-n_{k-q}}\geq\delta^{n_{k-q}}\}\cap\Omega_{k-q}$ is of
dimension at most equal to $q-1$ for some integer $n_{k-q}$.

When $q=0$ the last intersection is empty. We deduce that
$\kappa_{-n_k}<\delta_0^{n_k}$ on $\Omega_k$ for some constant
$\delta_0<\delta$. This, the fact that $f^{-1}(\Omega_k)\subset
\Omega_k$ together with the properties of cocycles imply
that $\kappa_{-mn_k}\leq \delta_0^{mn_k}$ on $\Omega_k$ for $m\geq
0$. By Lemma \ref{lemma_cocycle_3}, $\limsup
[\kappa_{-n}]^{1/n}\leq \delta_0$ on $\Omega_k$.

If $x$ is an arbitrary point in $X$, define 
$\kappa_-(x):=\liminf[\kappa_{-n}(x)]^{1/n}$. When
$\kappa_-(x)<\delta$,  $x$ belongs to 
 $\Omega_k$. Hence, $\limsup [\kappa_{-n}(x)]^{1/n}\leq
\delta_0<\delta$. This holds for any $\delta>\kappa_-(x)$. We
deduce that $[\kappa_{-n}(x)]^{1/n}$ converge to $\kappa_-(x)$.
We also obtain from the previous analysis that
$\{\kappa_-\geq\delta\}$ is equal to $\Sigma_1\cup\ldots\cup \Sigma_k$.
By Proposition \ref{prop_cocycle}, this set is contained
in the orbit of $\{\kappa_n\geq \delta^n\}$ for every $n\geq 0$.
\hfill $\square$

\

\noindent
{\bf Proof of Corollary \ref{cor_cocycle}.} Let $x$ be a point in $X$.
Assume that  $x_m:=f^{mn_k}(x)$ is in
$\Omega_k$ for every $m\geq 0$. We have $\kappa_{-n_k}(x_m)\leq
\delta_0^{n_k}$ for every $m\geq 0$ and hence
$\kappa_{n_k}(x_{m-1})\leq \delta_0^{n_k}$ for every $m\geq 1$.
Since $(\kappa_n)$ is sub-multiplicative, 
$\kappa_{mn_k}(x)\leq \delta_0^{mn_k}$ for every $m\geq 1$. Lemma
\ref{lemma_cocycle_bis} implies that
$\limsup[\kappa_n(x)]^{1/n}\leq\delta_0<\delta$. 
The first assertion follows.

Now assume that $(\kappa_n)$ is multiplicative. 
Define $\kappa_+(x):=\liminf [\kappa_n(x)]^{1/n}$. We have seen 
that if $x$ belongs to $\Sigma_1\cup\ldots\cup \Sigma_k$ then
$\kappa_+(x)\geq\delta$. So, the properties of multiplicative 
cocycles imply that if $x$ is sent by a
$f^n$ to a point in $\Sigma_1\cup\ldots\cup \Sigma_k$ then
$\kappa_+(x)\geq \delta$.
In other words, for any $\delta>\kappa_+(x)$ and $m\geq 0$, 
$x_m:=f^{mn_k}(x)$ is in $\Omega_k$. Therefore, we have 
$\limsup[\kappa_n(x)]^{1/n}\leq\delta_0<\delta$. This implies that 
$\limsup[\kappa_n(x)]^{1/n}\leq\kappa_+(x)$.
Hence, $[\kappa_n(x)]^{1/n}$ converge to $\kappa_+(x)$.
The identity $\kappa_+\circ f=\kappa_+$ is a consequence of the
properties of multiplicative cocycles.
\hfill $\square$

\begin{proposition} Let $(\kappa_n)$ be a sub-multiplicative cocycle
  as above. Let $x$ be a point in $X$. Then, for $l\geq 0$ large
  enough, $[\kappa_n(f^l(x))]^{1/n}$ converge to a constant
  $\kappa(x)$ which does not depend on $l$. In particular, $x$ is sent
  by an iterate of $f$ to the invariant analytic set $\{\kappa_-\geq
  \kappa(x)\}$. 
\end{proposition}
\proof
Let $\chi_l:=\limsup_{n\rightarrow
  \infty}[\kappa_n(f^l(x))]^{1/n}$. We have $\chi_l\geq 1$. Inequalities
(\ref{eq_cocycle}) for $s=1$ imply that $\kappa_{n+1}(f^l(x))\leq
M\kappa_n(f^{l+1}(x))$. Therefore, 
the sequence $(\chi_l)_{l\geq 0}$
is increasing. Define $\delta:=\lim_{l\rightarrow\infty}\chi_l=\sup_l\chi_l$. If
$\delta=1$, we have $\chi_l=1$; hence $[\kappa_n(f^l(x))]^{1/n}$
converge to 1 for every $l\geq 0$. Assume now that $\delta>1$. We use the
notations and the arguments in the proofs of Theorem
\ref{th_cocycle} and Corollary \ref{cor_cocycle}. When $l$ is large enough, we have
$\chi_l>\delta_0$. Therefore, $f^l(x)$ is sent by an iterate of $f$ to a point in
$\Sigma_1\cup\ldots\cup\Sigma_k$. So, for $l$ large enough, $f^l(x)$
belongs to $\Sigma_1\cup\ldots\cup\Sigma_k$. For such an $l$, we have
$\liminf_{n\rightarrow\infty}
[\kappa_n(f^l(x))]^{1/n}\geq\delta$. This implies that
$[\kappa_n(f^l(x))]^{1/n}$ converge to $\delta$. The last assertion in
the proposition is deduced from Corollary \ref{cor_cocycle}. 
\endproof

The following example shows that $[\kappa_n]^{1/n}$ do 
not always converge when $n$ goes to infinity.

\begin{example} \rm
Let $a$ be a fixed point of $f$  and $b$  a point in 
$f^{-1}(a)\setminus\{a\}$. 
Let $\lambda_n$ be positive numbers such that
$\lambda_{n+1}\leq\lambda_n+1$. 
Define the analytic sub-multiplicative cocycle
$(\kappa_n)$ by $\kappa_n(a):=e^n$, $\kappa_n(b):=e^{\lambda_n}$ and $\kappa_n(x):=1$
for $x\not=a,b$. 
Then, $[\kappa_n(b)]^{1/n}$ converge if and only if the arguments of
the complex numbers $n+i\lambda_n$
converge.
It is not difficult to find $\lambda_n$ which do not satisfy this
condition. In this case, the cocycle $(\kappa_n)$ is not
multiplicative. When $\lambda_n=0$, $[\kappa_n]^{1/n}$ converge but
the limit is not an invariant function.
\end{example}

\section{Exceptional sets for meromorphic maps} \label{section_dynamics}

There are some natural cocycles over a complex dynamical
system. We give some examples. Let
$f:\P^k\rightarrow\P^k$ be a holomorphic map of algebraic degree
$d\geq 2$ on the projective space of dimension $k$. It is induced by a map
$F:\C^{k+1}\rightarrow\C^{k+1}$ whose  irreducible components are homogeneous polynomials of
degree $d$. 
The map $f$ defines a ramified covering of degree $d^k$ on $\P^k$. 
We refer to the survey \cite{Sibony} for the basic
properties of $f$. 

Let $\kappa_n(x)$ denote the multiplicity of $f^n$ at $x$,
i.e. $f^{-n}(z)$ contains exactly $\kappa_n(x)$ points near $x$ for any $z$
generic close enough to $f^n(x)$. It is not difficult to check that
$(\kappa_n)$ is an analytic multiplicative cocycle. 
We now give an
application of Theorem \ref{th_cocycle}.

\begin{theorem}[\cite{DinhSibony1,DinhSibony3}] \label{th_except}
Let $f:\P^k\rightarrow\P^k$ be a holomorphic map of algebraic degree
$d\geq 2$. Then there is a proper analytic subset $\Ec$ of $\P^k$,
possibly empty, 
which is totally invariant by $f$:
$f^{-1}(\Ec)=f(\Ec)=\Ec$, and is maximal: $\Ec$ contains all the
proper analytic subsets $E$ of $\P^k$ satisfying $f^{-m}(E)\subset E$ for some
$m\geq 1$.
\end{theorem}
\proof
We use the multiplicative cocycle $(\kappa_n)$ defined above.
Let $E$ be as in the theorem. The sequence of analytic sets $(f^{-nm}(E))$ is
decreasing. Hence, there is an $n$ such that
$f^{-nm-m}(E)=f^{-nm}(E)$. Since $f^{nm}$ is surjective, we deduce
that $f^{-m}(E)=E$ and $f^m(E)=E$. Observe also
that if $\widetilde E:=E\cup\ldots\cup f^{-m+1}(E)$ then
$f^{-1}(\widetilde E)=\widetilde E$. Replacing $E$ by $\widetilde E$
allows to assume that $m=1$ and that $E$ is totally invariant: $f^{-1}(E)=f(E)=E$.  
The union of the irreducible components of $E$, of a given fixed dimension, is also totally
invariant. So, it is enough to consider the case where $E$ is of pure
dimension $p$. Let $n_0\geq 1$ be an integer such that $f^{n_0}$ fixes all the  irreducible components
of $E$. 

Denote by $g$ the restriction of
$f^{n_0}$ to $E'$, an  irreducible  component of $E$. We claim that the topological degree of $g$ is
equal to $d^{pn_0}$, that is, $g$ defines a ramified covering of degree
$d^{pn_0}$. Indeed, if $\omega$ is a K{\"a}hler form on $\P^k$, then
$(f^{n_0})^*(\omega^p)$ is cohomologous to $d^{pn_0}\omega^p$. By B{\'e}zout theorem, the
mass of the measure $g^*(\omega^p_{|E'})=(f^{n_0})^*(\omega^p)_{|E'}$
is equal to $d^{pn_0}$ times
the mass of $\omega^p_{|E'}$. 
This implies the claim.

Now, observe that the topological degree of $f^{n_0}$ is $d^{kn_0}$ while 
the topological degree of $g$ is equal to $d^{pn_0}$. This implies that
$\kappa_{-n_0}\geq d^{(k-p)n_0}$ on $E'$ since $f^{-n_0}(E')= E'$. The inequality still holds
if we replace $n_0$ by a multiple of $n_0$. Therefore, $E'$ is
contained in $\{\kappa_-\geq d^{k-p}\}$ and $E$ is contained in
$\Ec':=\{\kappa_-\geq d\}$. Since $E$ is totally invariant, it is
contained in $\Ec:=\cap_{n\geq 0} f^{-n}(\Ec')$. As above, we prove
that $\Ec$ is totally invariant. This completes the proof.
\endproof

The above result still holds for the restriction of $f$ to an
invariant irreducible 
analytic set $X$ of dimension $p$. One lifts $f$ to a
normalization of $X$. The map we obtain is a ramified covering of degree
$d^p$. The cocycle $(\kappa_n)$ over this map can be defined as
above. A simple decreasing induction on dimension allows to prove that $f$ admits a
finite number of totally invariant analytic subsets. Moreover, the
convex set of totally invariant probability measures, i.e. probability
measures $\nu$
such that $f^*(\nu)=d^k\nu$, is of finite dimension. The 
{\it equilibrium measure} $\mu$ is the only totally invariant
probability measure which has no mass on $\Ec$. The
exceptional set $\Ec$ is also characterized by the property that
$x\not\in\Ec$ if and only if $d^{-kn}(f^n)^*(\delta_x)$ 
converge to $\mu$ \cite{DinhSibony1, DinhSibony3}, see also
\cite{BriendDuval, FornaessSibony}. 
Here, $\delta_x$ denotes the Dirac mass at $x$.

Consider now a more general case. Let $(X,\omega)$ be a compact
K\"ahler manifold of dimension $k$ and let $f:X\rightarrow X$ be a dominant
meromorphic map. So, $f$ is a holomorphic map from $X\setminus I$ to $X$
where $I$ is the set of indeterminacy which is analytic and of
codimension $\geq 2$. The closure $\Gamma$ of the graph of
$f:X\setminus I\rightarrow X$ in $X\times X$ is an analytic set. If
$\pi_i:X\times X\rightarrow X$, $i=1,2$,  are the natural projections,
then $f$ is equal to $\pi_2\circ (\pi_{1|\Gamma})^{-1}$ on $X\setminus
I$. Define
for $A\subset X$
$$f(A):=\pi_2(\pi_1^{-1}(A)\cap\Gamma)\quad \mbox{and} \quad
f^{-1}(A):=\pi_1(\pi_2^{-1}(A)\cap\Gamma).$$
Define also $f^n:=f\circ \cdots\circ f$, $n$ times, on a suitable
dense Zariski
open set and extend it to a meromorphic map on $X$. In general, we do not
have $f^2(A)=f(f(A))$ nor $f^{-2}(A)=f^{-1}(f^{-1}(A))$. 

The set of indeterminacy $I$ is the set of points $x$ such that $f(x)$
is of dimension at least equal to $1$. 
Define $I_1:=I$ and $I_{n+1}:=f^{-1}(I_n)$ for $n\geq 1$. The orbit of
a point $x$
$$x,f(x), f^2(x),\ldots, f^n(x),\ldots$$
is well-defined if $x\not\in \cup_{n\geq 1} I_n$.

Let $I'$ denote the set of points $x$ such that $f^{-1}(x)$ is of
dimension larger than or equal to $1$. This is also an analytic subset of codimension
$\geq 2$ of $X$. The fibers 
$f^{-1}(x)$ for $x\not\in I'$, are finite and contain
the same number of points counted with multiplicity. This number $d_t$
is called {\it topological degree} of $f$. If $\nu$ is a finite measure on
$X$ having no mass on $I'$, then $f^*(\nu)$ is well-defined by the
formula
$$f^*(\nu):=(\pi_1)_*((\pi_{2|\Gamma})^*\nu).$$
If $\nu$ is positive, the mass of $f^*(\nu)$ is equal to $d_t$ times
the mass of $\nu$. Define $I'_1:=I'$ and $I'_{n+1}:=f(I'_n)$ for
$n\geq 1$. If
$x\not\in \cup_{n\geq 1} I_n'$, then $f^{-n}(x)$ contains exactly $d_t^n$
points counting with multiplicity. 

Choose a finite or countable union $\widetilde I$ of proper
analytic subsets of $X$ containing $I\cup I'$ such that 
if $\widetilde X:=X\setminus \widetilde I$ then 
$f^{-1}(\widetilde X)=f(\widetilde X)=\widetilde X$.
So,  $f:\widetilde X\rightarrow\widetilde X$ defines a ``ramified
covering'' of degree $d_t$.
We consider 
the topology on $\widetilde X$ induced by 
the Zariski topology on
$X$. As in the case of holomorphic maps on $X$, one can define cocycles
over $f:\widetilde X\rightarrow \widetilde X$.
Theorem \ref{th_cocycle} and Corollary \ref{cor_cocycle} still hold in
this case. 

Recall that the {\it dynamical degree} of order $p$ of $f$, $0\leq p\leq k$,
is defined by
$$d_p:=\lim_{n\rightarrow\infty} \Big(\int_X
(f^n)^*(\omega^p)\wedge\omega^{k-p}\Big)^{1/n}.$$
Note that $(f^n)^*(\omega^p)$ is a $(p,p)$-form with $L^1$
coefficients. The last limit always exists and the dynamical degrees
are bi-meromorphic invariants \cite{DinhSibony2}. We also have $d_0=1$
and $d_k=d_t$, the topological degree of $f$.

\begin{theorem} \label{th_excep_mer}
Let $f:X\rightarrow X$ be a dominant meromorphic map 
as above. Assume that the topological degree
of $f$ is strictly larger than the other dynamical degrees. Then there
is a proper analytic subset $\Ec$ of $X$, possibly empty, such that
\begin{enumerate}
\item all the irreducible components of $\Ec$ intersect $\widetilde X$;
\item $\Ec$ is totally invariant: $f^{-1}(\Ec\cap \widetilde X)=f(\Ec\cap
      \widetilde X)=\Ec\cap \widetilde X$;
\item $\Ec$ is maximal: if $E\subset X$ is a proper analytic subset  whose irreducible
      components intersect $\widetilde X$ and if $f^{-m}(E\cap
      \widetilde X)\subset E$ for some $m\geq 1$, then
      $E\subset\Ec$. 
\end{enumerate}
\end{theorem}

One can also characterize $\Ec$ by the following property. For $x$
a point in $\widetilde X$, the sequence of probability measures 
$d_t^{-n}(f^n)^*(\delta_x)$ converges to the
equilibrium measure of $f$ if and only if $x\not\in\Ec$, see
\cite{Dinh, DinhSibony3, Guedj} for the proof. 

The proof of Theorem \ref{th_excep_mer} follows the arguments in Theorem
\ref{th_except}. We only need the following lemma.

\begin{lemma} Let $E$ be a proper irreducible analytic subset of $X$ which
  intersects $\widetilde X$ and $n_0\geq 1$ an integer such that
  $f^{n_0}(E\cap \widetilde X)\subset E$. Then the topological degree of
  the meromorphic map $f^{n_0}:E\rightarrow E$ is at most equal to
  $\delta^{n_0}$ where $0<\delta<d_t$ is a constant independent of $E$
  and $n_0$. 
\end{lemma}
\proof
Denote by $g$ the map $f^{n_0}:E\rightarrow E$.
Let $p$ be the dimension of $E$, $d$ the topological degree of
$g$ and $\delta<d_t$ a constant larger
than or equal to the dynamical degrees $d_0,\ldots, d_{k-1}$. The mass
of $g^*(\omega^p_{|E})$ is equal to $d$ times the mass of
$\omega^p_{|E}$. 
We will bound the mass of $g^*(\omega^p_{|E})$. Fix a constant $\epsilon>0$.
In what follows, all the constants are independent of $E$, $n_0$ and $n$.

The form $(f^{n_0})^*(\omega^p)$ has $L^1$ coefficients and is 
smooth on the Zariski open set $X\setminus \cup_{1\leq n\leq n_0} I_n$
which contains $\widetilde X$. 
It also defines a positive closed $(p,p)$-current on $X$ whose mass is
given by the formula
$$\|(f^{n_0})^*(\omega^p)\|:=\int_X (f^{n_0})^*(\omega^p)\wedge
\omega^{k-p}.$$
Therefore, there is a constant $c>0$ 
such that this mass is bounded by $c(d_p+\epsilon)^{n_0}$. We
regularize $(f^{n_0})^*(\omega^p)$ using the results in
\cite{DinhSibony2}. There are smooth positive closed $(p,p)$-forms
$\Omega_i^\pm$ such that $\Omega^+_i-\Omega^-_i$ converge to $(f^{n_0})^*(\omega^p)$
in the sense of currents and 
$$\|\Omega_i^\pm\|\leq A\|(f^{n_0})^*(\omega^p)\|\leq
A'(d_p+\epsilon)^{n_0}$$
for some constants $A$ and $A'$. Since $(f^{n_0})^*(\omega^p)$ is smooth
in a neighbourhood of $\widetilde X$, the construction in
\cite{DinhSibony2} gives forms $\Omega_i^\pm$ such that
$\Omega_i^+-\Omega_i^-$ converge uniformly on compact subsets of
$\widetilde X$ to $(f^{n_0})^*(\omega^p)$. It follows that
$$\|g^*(\omega^p_{|E})\|=\|(f^{n_0})^*(\omega^p)_{|E\cap\widetilde X}\|\leq \limsup_{i\rightarrow\infty} \int_E
\Omega_i^+.$$
The last integral depends only on the cohomology class of
$\Omega_i^+$. The above bound on the mass of $\Omega_i^+$ induces a bound
on its cohomology class. Therefore, $\|g^*(\omega^p_{|E})\|\leq
A''(d_p+\epsilon)^{n_0}$ for some constant $A''$. 
The estimate also holds for the iterates of $g$: we have 
$\|(g^n)^*(\omega^p_{|E})\|\leq
A''(d_p+\epsilon)^{nn_0}$. It follows that $d^n\lesssim
(d_p+\epsilon)^{nn_0}$. Hence, $d$ is at most equal to
$d_p^{n_0}$. This completes the proof.
\endproof


\noindent
T.-C. Dinh, UPMC Univ Paris 06, UMR 7586, Institut de
Math{\'e}matiques de Jussieu, F-75005 Paris, France. {\tt
  dinh@math.jussieu.fr}, {\tt http://www.math.jussieu.fr/$\sim$dinh} 

\end{document}